\documentclass[reqno]{amsart}
\usepackage{amssymb,amscd,enumerate,mathrsfs,curves}

\numberwithin{equation}{section}
\newtheorem{theorem}[equation]{Theorem}

\newtheorem{lemma}[equation]{Lemma}

\theoremstyle{definition}

\newtheorem{example}[equation]{Example}
\DeclareMathOperator{\GAMDiff}{Diff}
\DeclareMathOperator{\GAMdist}{dist}
\DeclareMathOperator{\GAMev}{ev}
\DeclareMathOperator{\GAMHom}{Hom}
\DeclareMathOperator{\GAMrank}{rank}
\DeclareMathOperator{\GAMspec}{spec}
\DeclareMathOperator{\GAMsym}{\pmb\sigma}

\def\GAMagen{\mathfrak a}
\def\GAMbgen{\mathfrak b}
\def\GAMExt{\mathfrak E}
\def\GAMHol{\mathfrak{H}}
\def\GAMKK{\mathfrak{K}}
\def\GAMm{\mathfrak m}
\def\GAMMero{\mathfrak{M}}
\def\preP{\mathfrak P}
\def\ssp{\mathfrak s}

\def\GAMC{\mathbb C}
\def\GAMNN{\mathbb N}
\def\GAMR{\mathbb R}

\def\GAMDom{\mathcal D}
\def\GAMM{\mathcal M}
\def\GAMN{\mathcal N}
\def\GAMRing{\mathcal R}
\def\GAMY{\mathcal Y}
\def\GAMZ{\mathcal Z}

\def\GAMtarget{\mathscr G}
\def\GAMkerb{\mathscr K}
\def\GAMtrt{\mathscr S}
\def\GAMtrb{\mathscr T}

\def\GAMminus{\backslash}
\def\GAMim{i}
\def\GAMId{I}
\def\GAMembed{\hookrightarrow}

\def\GAMopen#1{\smash[t]{\overset{{}_{\,\,\,\circ}}{#1}{}}}
\def\GAMset#1{\{#1\}}
\def\GAMjb#1{\langle#1\rangle}

\def\GAMrpar{)}
\def\GAMlbra{[}

\def\GAMwC{\,{}^w{\!}C}
\def\GAMwev{\,{}^w\hspace{-0.5pt}{\GAMev}}
\def\GAMwT{\,{}^w\hspace{-0.5pt}T}
\def\GAMwev{\,{}^w\hspace{-0.5pt}{\GAMev}}
\def\GAMwpi{\,{}^w\!\pi}
\def\GAMwsym{\,{}^w\!\!\GAMsym}

\def\GAMeT{\,{}^e\hspace{-0.5pt}T}

\def\GAMchern{\mathrm{c}}

\def\GAMWedge{\raise2ex\hbox{$\mathchar"0356$}}
\def\GAMbT{\,{}^b\hspace{-0.5pt}T}
\def\GAMbP{\,{}^b\!P}
\def\GAMbPhat{\,{}^b\!\widehat P}
\def\GAMbPtilde{\,{}^b\!\widetilde P}
\def\GAMA{\,{}^b\!A}
\def\GAMie{i.~e.}
\def\GAMopcit{op.~cit.}
\def\GAMbysame{\underline{\hspace*{20pt}}}

\begin{document}

\title[Boundary value problems]{Boundary value problems for elliptic wedge operators: the first order case}
\thanks{Work partially supported by the National Science Foundation, Grants DMS-0901202 and DMS-0901173}

\author{Thomas Krainer}
\address{Penn State Altoona\\ 3000 Ivyside Park \\ Altoona, PA 16601-3760}
\email{krainer@psu.edu}
\author{Gerardo A. Mendoza}
\address{Department of Mathematics\\ Temple University\\ Philadelphia, PA 19122}
\email{gmendoza@temple.edu}

\begin{abstract}
This note is a description of some of the results obtained by the authors in connection with the problem in the title. These, discussed following a summary  of background material concerning wedge differential operators, consist of the notion of trace bundle, an extension of the Douglis-Nirenberg calculus to handle spaces of anisotropic varying regularity and associated pseudodifferential operators, and boundary value problems proper, the latter in the first order case. The concepts concerning the main results are illustrated with simple examples.

\end{abstract}



\maketitle

%%%%%%%%%%%%%%%%%%%%%%%%%%%%%%%%%%%%%%%%%%%%%%%%%
%%%%%%%%%%%%%%%%%%%%%%%%%%%%%%%%%%%%%%%%%%%%%%%%%

\section{Introduction}%\label{sec-introduction}

The present note is an account of results published in a series of papers \cite{GAMKrMe12a,GAMKrMe12b,GAMKrMe13} in connection with boundary value problems for elliptic wedge operators on a manifold with fibered boundary. Briefly, in \cite{GAMKrMe12a} we  address the fundamental issue of boundary values, in \cite{GAMKrMe12b} we construct an extension of the Douglis-Nirenberg calculus (see for instance \cite{GAMChazarainPiriou,GAMHor67} for the role of this calculus in the classical context), while in \cite{GAMKrMe13} we address elliptic boundary value problems for first order wedge operators and prove, in particular, sufficient conditions for well-posedness of such problems. Here we shall address the main aspects of each of these papers in subsequent sections, after dealing with background information and some notation. 

Our point of view, properly translated, closely parallels that of regular elliptic boundary value problems, and indeed our approach, restricted in \cite{GAMKrMe13} to first order operators, allows a full analysis of the classical problem as a special case.

Of course the work described here does not exist in a vacuum. However, we shall cite only the most directly pertinent work, and refer the reader to the papers indicated above for a more representative listing of research in the general area of elliptic problems for elliptic operators on manifolds with singularities. The preceding notwithstanding, we call attention to the paper of Mazzeo and Vertman \cite{GAMMazzVertCommunication} treating higher order problems under certain assumptions.

%\begin{acknowledgement}
%Work partially supported by the National Science Foundation, Grants DMS-0901202 (TK) and DMS-0901173 (GAM).
%\end{acknowledgement}

%%%%%%%%%%%%%%%%%%%%%%%%%%%%%%%%%%%%%%%%%%%%%%%%%
\section{Set-up}\label{sec-setup}

The differential-topological setup is that of Mazzeo \cite{GAMMazz91}. Namely, a compact manifold $\GAMM$ whose boundary $\GAMN$ is the total space of a locally trivial fibration $\wp:\GAMN\to\GAMY$ with typical fiber $\GAMZ$. (For an extension of this kind of structure the reader is directed to the paper by Albin, Leichtnam, Mazzeo, and Piazza, \cite{GAMAlbinLeichtMazzPiazza12} which gives a very clear description of the process of resolving singularities of an arbitrary stratified pseudomanifold through a series of blowups keeping track of boundary fibrations.) The base space $\GAMY$, the edge, may have several components which are manifolds of possibly different dimensions and the typical fibers over different components may not be diffeomorphic, but we will ignore this for notational simplicity.

The analytic objects are wedge operators, \GAMie, elements of $x^{-m}\GAMDiff^m_e(\GAMM;E,F),$ where $\GAMDiff^m_e(\GAMM;E,F)$ is the class of edge differential operators of order $m$ defined by Mazzeo, that is, linear differential operators on $\GAMM$ of order $m$ with smooth coefficients which along the boundary differentiate only in directions tangent to the fibers; $E$ and $F$ are Hermitian vector bundles over $\GAMM$, and  $x$ is a defining function for $\GAMN$, positive in $\GAMopen \GAMM$. An element of $\GAMDiff^m_e(\GAMM;E,F)$ is thus a regular differential operator on $\GAMM$ which in local coordinates $x,y_j,z_\mu$ near any point of the boundary, with $x$ as just indicated and the $y_j$ restricted to $\partial\GAMM$ being constant of fibers, has the form
\begin{equation}\label{LocalP}
\sum_{k+|\alpha|+|\beta|\leq m}a_{k\alpha\beta}(xD_x)^k (xD_y)^\alpha D_z^\beta
\end{equation}
with respect to local trivializations of $E$ and $F$. The coefficients  $a_{k\alpha\beta}$ are smooth up to the boundary. A regular differential operator $A$ of order $m$ on $\GAMM$ is an element of $x^{-m}\GAMDiff^m_e$ by way of the cheap trick $A=x^{-m}(x^mA)$. In this case $\GAMY=\GAMN$ and the fibers are just the points of $\GAMN$. One of the initial motivations for the structural specification of the elements of $\GAMDiff^m_e$ comes from what results when writing a regular differential operator in cylindrical coordinates.

The functional analytic component enters through a choice of a $b$-density $\GAMm_b=x^{-1}\GAMm$ as in Melrose \cite{GAMMel93}; $\GAMm$ is a smooth density. With the $b$-density and the Hermitian structures of $E$ and $F$ one gets weighted $L^2$ spaces, e.g. $x^{-\gamma}L^2_b(\GAMM;E)$.

%%%%%%%%%%%%%%%%%%%%%%%%%%%%%%%%%%%%%%%%%%%%%%%%%
\section{Some considerations}\label{sec-considerations}

Any elliptic element of $\GAMDiff^m_e(\GAMM;E,F)$ (we review the intrinsic notion of ellipticity in Section~\ref{sec-ellipticity}), viewed initially as an operator
\begin{equation}\label{Core}
C_c^\infty(\GAMopen\GAMM;E)\subset x^{-\gamma}L^2_b(\GAMM;E)\to x^{-\gamma}L^2_b(\GAMM;F),
\end{equation}
admits only one closed extension ($\gamma$ is a real number; the spaces are $L^2$ spaces with respect to the measure $x^{2\gamma} \GAMm_b$), whereas generically elliptic elements of the space $x^{-m}\GAMDiff^m_e(\GAMM;E,F)$ admit infinitely many such extensions: this is the reason why boundary value problems make sense for wedge operators but not for edge operators. Having a unique closed extension is a property shared by other classes of operators such as the differential operators in the $\Theta$-calculus of Epstein, Melrose, and Mendoza \cite{GAMEpRBMMe91}, and more generally, those associated to Lie structures at infinity of Ammann, Lauter, and Nistor \cite{GAMAmLaNi07}.

Let $A$ be an elliptic element in $x^{-m}\GAMDiff^m_e(\GAMM;E,F)$. Recall that the domain of the maximal extension of $A$, initially as an operator \eqref{Core}, is the space 
\begin{equation*}
\GAMDom_{\max}(A)=\GAMset{u \in x^{-\gamma}L^2_b(\GAMM;E): Au\in x^{-\gamma}L^2_b(\GAMM;F)}
\end{equation*}
and that the minimal domain, $\GAMDom_{\min}(A)$, is the domain of the closure of $A$ starting with \eqref{Core}. In the case of a regular elliptic operator of order $m$, the minimal domain is $H^m_0(\GAMM;E)$ (we take $\gamma=1/2$ in this case because $L^2(\GAMM;\GAMm)=x^{-1/2}L^2(\GAMM,\GAMm_b)$).

One seeks among other things to establish the existence of a split exact sequence
\begin{equation}\label{SplitExactSequence}
0 \to \GAMDom_{\min}(A) \to H_A^m \to {\mathcal S}_A \to 0
\end{equation}
in which $H^m_A$ is a conveniently chosen subspace of the maximal domain. For regular elliptic differential operators this sequence is analogous, and in a natural general theory should reduce, to the classical sequence
\begin{equation*}
0 \to H^m_0(\GAMM,E) \to H^m(\GAMM,E) \to \bigoplus_{j=0}^{m-1}H^{m-j-1/2}(\partial\GAMM;E_{\partial\GAMM}) \to 0
\end{equation*}
that is associated with taking Cauchy data on the boundary. On the face of it, in \eqref{SplitExactSequence} one could take $H_A^m=\GAMDom_{\max}(A)$ and ${\mathcal S}_A=\GAMDom_{\max}(A)/\GAMDom_{\min}(A)$; viewing the quotient as the orthogonal of $\GAMDom_{\min}(A)$ in $\GAMDom_{\max}(A)$ with respect to the inner product 
\begin{equation}\label{InnerProd}
(u,v)_A=(u,v)+(Au,Av)
\end{equation}
on the maximal domain gives a natural splitting. However this choice is generally bad, because one should also require that the inclusion $H_A^m\GAMembed x^{-\gamma}L^2_b$  be compact. The following example shows that this need not be the case when $H^m_A$ is taken to be $\GAMDom_{\max}(A)$ endowed with the norm defined by the inner product \eqref{InnerProd}.

\begin{example}\label{BasicExample}
Let $\GAMM$ be the closed unit disk in $\GAMR^2$ and let $\Delta$ be the Euclidean Laplacian. We claim that the inclusion of $\GAMDom_{\max}(\Delta)$ in $L^2(\GAMM)$ is not compact. If it were, then also the inclusion of $\ker\Delta$ in $L^2$ is compact. But with the norm defined by \eqref{InnerProd} in $\GAMDom_{\max}(\Delta)$ we have
\begin{equation}\label{DeltaNorm}
\|u\|_{\Delta}^2=\|u\|^2+\|\Delta u\|^2,
\end{equation}
so the $\Delta$-norm and the $L^2$ norm are the same on $\ker \Delta$. But the compactness of the inclusion map now implies that the unit sphere of $\ker\Delta$ is compact in the $L^2$ norm, a contradiction since $\ker\Delta$ is infinite-dimensional. Thus the inclusion of $\GAMDom_{\max}(\Delta)$ in $L^2(\GAMM)$ is not compact.
\end{example}

We now discuss briefly the role of the weight $x^{2\gamma}$ and the use of a $b$-density rather than a regular density. The weight $x^{2\gamma}$ connects with geometric information such as what appears when introducing cylindrical coordinates along a submanifold $\GAMY$ of codimension $k$ in a smooth manifold: a smooth measure near $\GAMY$ becomes essentially $x^{k-1}dx\,dy\,dz$, $x$ being the radial variable and $dz$ representing the measure on the sphere $S^{k-1}$. Both the factor $x^{-1}$ making up the $b$-density and the weight $x^{2\gamma}$ can be removed by conjugating the operator with multiplication by an appropriate power of $x$. We eventually take advantage of this and pick $\gamma=m/2$, but keep the $b$-density since this brings to the foreground the multiplicative structure of $\GAMR_+$ (for which $x^{-1}dx$ is a Haar measure). The class $x^{-m}\GAMDiff^m_e$ is invariant under conjugation as described, however not so the class $\GAMDiff^m$ of regular differential operators, which under such operations end up subsumed in the more general class of wedge operators; depending on the particularities of the problem, this may be advantageous.

%%%%%%%%%%%%%%%%%%%%%%%%%%%%%%%%%%%%%%%%%%%%%%%%%
\section{Ellipticity, the wedge cotangent bundle and the structure ring}\label{sec-ellipticity}

Ellipticity of an element $A\in x^{-m}\GAMDiff^m_e(\GAMM;E,F)$ means ellipticity of $P=x^m A$ in the sense of \cite{GAMMazz91}. This is ellipticity of $P$ over the interior of $\GAMM$ in the usual sense, and, near a boundary point where $P$ is written as in \eqref{LocalP}, invertibility of 
\begin{equation*}
\sum_{k+|\alpha|+|\beta|= m}a_{k\alpha\beta}\xi^k \eta^\alpha \zeta^\beta.
\end{equation*}
While this is a perfectly good practical definition of ellipticity, it disregards basic information of the manifold-with-boundary-fibration and the class of $w$-differential operators. Still, fixing $x$ allows a definition of principal symbol of $A$ by way of declaring it to be the edge-principal symbol of $x^mA$, which is what the above expression is in coordinates. This is, however, not quite satisfactory since it does depend, albeit mildly, on the choice of defining function.

The natural structure bundle for the class of wedge differential operators is the 
wedge cotangent bundle $\GAMwT^*\GAMM$. It is constructed in a fashion similar to Melrose's $b$-tangent bundle in \cite{GAMMel81,GAMMel93} or Mazzeo's edge tangent bundle in \cite{GAMMazz91}, as follows (see \cite{GAMGiKrMe10} for details). The space of continuous differential $1$-forms on $\GAMM$ whose pull-back to the fibers of $\GAMN\to\GAMY$ vanishes is a finitely generated projective module, $\GAMwC(\GAMM;T^*\GAMM)$, over the ring of continuous functions on $\GAMM$, and is therefore, by a theorem of Swan \cite{GAMSwan62}, (isomorphic to) the space of sections of a vector bundle over $\GAMM$ which we denote by $\GAMwT^*\GAMM$. This vector bundle is easily seen to be a $C^\infty$ bundle, and is the natural structure bundle to the same extent as $\GAMeT\GAMM$ is the natural structure bundle in the case of edge operators and $\GAMbT\GAMM$ is in the case of $b$-operators. The inclusion map $\GAMwev_*:\GAMwC(\GAMM;T^*\GAMM)\to C(\GAMM;T^*\GAMM)$ determines a (smooth) bundle homomorphism $\GAMwev:\GAMwT^*\GAMM\to T^*\GAMM$ covering the identity; $\GAMwev$ is an isomorphism over the interior while over the boundary its kernel is the conormal bundle to the fibers of $\GAMZ\to\GAMY$. 

The naturality of $\GAMwT^*\GAMM$ is further justified by it being the domain of a principal symbol map for elements $A\in x^{-m}\GAMDiff^m_e(\GAMM;E,F)$: there is a smooth homomorphism
\begin{equation*}
\GAMwsym(A)\in C^\infty(\GAMwT^*\GAMM\GAMminus 0;\GAMHom(\GAMwpi^*E,\GAMwpi^*F)),
\end{equation*}
the wedge symbol of $A$, related to the standard principal symbol of $A$ over $\GAMopen \GAMM$ by 
\begin{equation*}
\GAMwsym(A) = \GAMsym(A)\circ\GAMwev.
\end{equation*}
Naturally, ellipticity is defined as invertibility of $\GAMwsym(A)$.

The section $\GAMwsym(A)$ can also be obtained by an oscillatory test using real-valued functions in the ring
\begin{equation*}
\GAMRing=\GAMset{f\in C^\infty(\GAMM):f|_\GAMN\text{ is constant on the fibers of }\wp},
\end{equation*}
see  \cite{GAMKrMe13}. The fundamental role of this ring can be seen from the observation that it determines the boundary structure of $\GAMM$. It can also be used to define the spaces $x^{-m}\GAMDiff^m_e$ without resorting to coordinates, as described in the introduction of the just cited paper. From another point of view, the differentials of real-valued elements of $\GAMRing$ generate $\GAMwC^\infty(\GAMM;T^*\GAMM)$ as a module over $C^\infty(\GAMM;\GAMR)$. Finally, observe that if the configuration of $\GAMM$ with its boundary fibration comes from blowing up a smooth manifold along a smooth submanifold (that is, from cylindrical coordinates), then $\GAMRing$ is, to first order along $\GAMN$, the pull-back of the ring of smooth functions on the original manifold.

Incidentally, when a regular elliptic operator on a smooth manifold is written in cylindrical coordinates with axis a given submanifold, the result is not just a wedge operator as already pointed out, but it is also a $w$-elliptic operator. Along the same vein, a smooth Riemannian metric on the cotangent bundle of the original manifold becomes a smooth metric on the wedge cotangent bundle.

%%%%%%%%%%%%%%%%%%%%%%%%%%%%%%%%%%%%%%%%%%%%%%%%%
\section{Indicial and normal families}%\label{sec-othersymbols}

The operator $P=x^mA$ is in particular a $b$-operator (by way of replacing the fine structure of the fibration $\GAMN\to\GAMY$ with the one in which each connected component of $\GAMN$ is one fiber): $P\in \GAMDiff^m_b(\GAMM;E,F)$. The operator $P$ has the property that if $\phi\in C^\infty(\GAMM;E)$ then $(P\phi)|_{\GAMN}$ depends only on $\phi|_\GAMN$, thus giving an operator $P\big|_\GAMN:C^\infty(\GAMN;E_\GAMN)\to C^\infty(\GAMN;F_\GAMN)$. By $E_\GAMN$ we mean the part of $E$ over $\GAMN$. Also $x^{-\GAMim\sigma}Px^{-\GAMim\sigma}\in \GAMDiff^m_b$, and the indicial family of $A$ (or $P$) is defined as 
\begin{equation*}
\GAMbPhat(\sigma)= (x^{-\GAMim\sigma}Px^{\GAMim\sigma})\big|_\GAMN, \quad\sigma\in \GAMC.
\end{equation*}
Because of the factor $x$ occurring with each derivative in $y$, this operator does not differentiate in $y$: it gives a family of differential operators $\widehat P_y(\sigma):C^\infty(\GAMZ_y;E_{\GAMZ_y})\to C^\infty(\GAMZ_y;F_{\GAMZ_y})$ depending smoothly on $(\sigma,y)\in \GAMC\times \GAMY$ and holomorphically (polynomially) in $\sigma$. Here $\GAMZ_y$ is the fiber of $\wp:\GAMN\to\GAMY$ over $y$.

Let $\pi_\wedge:\GAMN^\wedge\to \GAMN$ be the closed inward-pointing part of the normal bundle $N_\GAMM\GAMN$ of $\GAMN$ in $\GAMM$. Its boundary, being the zero section of $N_\GAMM\GAMN$, is identified with $\GAMN$ and so carries the same fibration structure as $\partial\GAMM$. The part of $\GAMN^\wedge$ over $\GAMZ_y$ is denoted $\GAMZ_y^\wedge$.

The infinitesimal generator of the radial action on $\GAMN^\wedge$ is $x_\wedge\partial_{x_\wedge}$ where $x_\wedge$ is any linear function $N_\GAMM\GAMN$, positive on $\GAMopen \GAMN^\wedge$, for example the one induced by $dx$. We shall write also $x$ for $x_\wedge$ since there is little risk of confusion. The family $\GAMbPhat(\sigma)$ is polynomial in $\sigma$, so replacing $\sigma$ by $-\GAMim x\partial_x$ (more properly, by $-\GAMim \nabla_{x\partial_x}$) gives an operator $\GAMbP\in \GAMDiff^m_b(\GAMN^\wedge;E^\wedge,F^\wedge)$ (which does not differentiate in $y$) where for example $E^\wedge$ is the pullback of $E$ to $\GAMN^\wedge$ (and $\nabla$ the pullback of some (fixed) Hermitian connection on $E$). The indicial operator of $A$ is then defined to be the operator
\begin{equation*}
\GAMA=\frac{1}{x^m}\GAMbP\in \frac{1}{x^m}\GAMDiff^m_b(\GAMN^\wedge;E^\wedge,F^\wedge).
\end{equation*}
In local coordinates, if $x^m A$ is given by \eqref{LocalP} near some boundary point, then 
\begin{equation*}
\GAMA=\frac{1}{x^m}\sum_{k+|\beta|\leq m}a_{k0\beta}(0,y,z)(xD_x)^kD_z^\beta.
\end{equation*}
Note that because $\GAMbP$ does not differentiate in $y$, the indicial operator may be viewed as a family
\begin{equation*}
\GAMY\ni y\mapsto \GAMA_y=\frac{1}{x^m}\GAMbP\in \frac{1}{x^m}\GAMDiff^m_b(\GAMZ_y^\wedge;E_{\GAMZ_y^\wedge},F_{\GAMZ_y^\wedge}).
\end{equation*}
We wrote $E_{\GAMZ^\wedge_y}$ to mean $E^\wedge$ along $\GAMZ_y^\wedge$ with a slight abuse of the notation. 

One defines an $\GAMR_+$-action $\varrho\mapsto \kappa_\rho$ on sections $u$ of $E^\wedge$ by letting $(\kappa_\varrho u)(\nu)$ be the result of parallel transport of $u(\varrho \nu)$ from $\varrho\nu$ to $\nu$ along the fiber of $\GAMN^\wedge$ followed by multiplication by $\varrho^\gamma$ (\GAMie, on functions, $(\kappa_\varrho f)(\nu)=\varrho^\gamma f(\varrho \nu)$). Thus $\kappa_\varrho$ is a bundle homomorphism covering the radial action $\nu\mapsto \tau_\varrho\nu =\varrho\nu$. The factor $\varrho^\gamma$ ensures that $\kappa_\varrho$ is unitary on $x^{-\gamma} L^2_b(\GAMN^\wedge;E^\wedge)$. The $b$-density for the latter space is defined with the aid of a tubular neighborhood map
\begin{equation*}
\GAMN^\wedge\supset V \xrightarrow{\varphi} W\subset \GAMM
\end{equation*}
as
\begin{equation*}
\GAMm_b^\wedge=\lim_{\varrho\to 0} (\varphi\circ\tau_\varrho)^*\GAMm_b
\end{equation*}

The normal family of $A$ associates to each element $\pmb\eta\in T^*Y\GAMminus 0$ an element of $x^{-m}\GAMDiff^m_b(\GAMZ_y^\wedge;E_{\GAMZ_y^\wedge},F_{\GAMZ_y^\wedge})$ by means of the formula 
\begin{equation}\label{NormalFamilyAsLimit}
A_\wedge (\pmb\eta)u= (\lim_{\varrho\to \infty} \varrho^{-m}\kappa^{-1}_\varrho e^{-\GAMim \varrho \wp_\wedge^*g} \Phi^* A \Phi_* e^{\GAMim \varrho \wp_\wedge^*g} \kappa_\varrho \tilde u)\big|_{\GAMZ_y^\wedge},\ u\in C_c^\infty(\GAMopen\GAMZ_y^\wedge;E_{\GAMZ_y^\wedge}),
\end{equation}
with $g\in C^\infty(\GAMY)$ real-valued with $dg(y)=\pmb\eta$, and $\tilde u$ a $C_c^\infty$ extension of $u$. The maps $\Phi_*$ (and its inverse $\Phi^*$) are defined using $\varphi$, the radial action on $\GAMN^\wedge$, and parallel transport on $E$ and $F$ as needed. See \cite[Proposition 2.10]{GAMGiKrMe10} for details on this; a local argument shows that $A_\wedge$ is independent of $\varphi$. In local coordinates, if $x^m A$ is given by \eqref{LocalP} near some boundary point, then 
\begin{equation*}
A_\wedge(\pmb\eta)=\frac{1}{x^m}\sum_{k+|\alpha|+|\beta|\leq m}a_{k\alpha\beta}(0,y,z)(xD_x)^k(x\eta)^\alpha D_z^\beta,\quad\pmb\eta=\eta\cdot dy.
\end{equation*}
This expression makes sense also at $\pmb\eta=0$, were it becomes canonically equal to $\GAMA$.

%%%%%%%%%%%%%%%%%%%%%%%%%%%%%%%%%%%%%%%%%%%%%%%%%
\section{The kernel bundle of a holomorphic Fredholm family}%\label{sec-kernelbundle}

The argument leading to the definition of $\GAMbP$ can be extended to give a formal power series expansion 
\begin{equation*}
P=\sum_{\ell=0}^N  \GAMbP_{\ell}x^\ell+\GAMbPtilde_{N+1}x^{N+1}
\end{equation*}
in which the $\GAMbP_{\ell}$ and $\GAMbPtilde_{N+1}$ are all $b$-operators on $\GAMN^\wedge$, the $\GAMbP_\ell$  commute with $\nabla_{x\partial_x}$ and $\GAMbPtilde_{N+1}$ is defined only near $\GAMN$. The operator $\GAMbP_{0}$ is equal to $\GAMbP$, each  operator $\GAMbP_{\ell}$ includes derivatives in $y$ of at most order $\max\GAMset{\ell,m}$, similarly the remainder $\GAMbPtilde_{N+1}$.

If $u\in x^{-\gamma}L^2_b(\GAMM;E)$, its Mellin transform is 
\begin{equation*}
\widehat u(\sigma,p)=\int_{\pi^{-1}_\wedge(p)} x^{-\GAMim\sigma}\omega(x)\Phi^*u \,\frac{dx}{x}
\end{equation*}
This may be viewed variously as an element of $L^2$ on $\GAMset{\Im \sigma=\gamma} \times \GAMN$, as a holomorphic function on $\Im\sigma>\gamma$ with values in $L^2(\GAMN,E|_\GAMN)$, and so on. Changing the cut-off function $\omega$ (which is supported in a small neighborhood of $\GAMN$) changes $\widehat u(\sigma,p)$ by an entire additive term.

One gets a hold on boundary values of elements $u\in\GAMDom_{\max}(A)$ by exploiting an idea from \cite{GAMMM83} here applied to the fact that for such $u$ one has that $\widehat u(\sigma,p)$ is holomorphic in $\Im\sigma>\gamma$ whereas $(x^mAu)\widehat{\ }=\widehat f(\sigma)$ is holomorphic in $\Im\sigma>\gamma-m$, so 
\begin{equation*}
\GAMbPhat (\sigma)\widehat u(\sigma)=\widehat f(\sigma)-\sum_{\ell=1}^N\GAMbPhat_\ell(\sigma)\widehat u(\sigma+\GAMim\ell)-(\GAMbPtilde_{N+1}x^{n+1} u)\widehat{\ }(\sigma),
\end{equation*}
where the right hand side is holomorphic in $\Im\sigma>\gamma-1$ (but in principle less regular in $y$ than $u$ since it includes derivatives in $y$). It follows that 
\begin{equation}\label{BeginExpansion}
\widehat u(\sigma)=\GAMbPhat (\sigma)^{-1}[\widehat f(\sigma)-\sum_{\ell=1}^N\GAMbPhat_\ell(\sigma)\widehat u(\sigma+\GAMim\ell)-(\GAMbPtilde_{N+1}x^{n+1} u)\widehat{\ }(\sigma)],
\end{equation}
for $\Im\sigma>\gamma-1$ such that $\GAMbPhat (\sigma)^{-1}$ exists. 

Assume henceforth that $A$ is elliptic. Then $\widehat P_y(\sigma)$ is elliptic for each $(\sigma,y)$ and invertible in regions $|\Im\sigma|<a$ for any $a$ when $|\Re \sigma|$ large enough, uniformly for $y$ in compact sets (so on $\GAMY$ itself). This implies in particular that the set 
\begin{equation*}
\GAMspec_{b,y}(A)=\GAMset{\sigma:\widehat P_y(\sigma)\text{ is not invertible}},
\end{equation*}
the boundary spectrum of $A$ (or $P)$ at $y$ is discrete (see \cite{GAMMM83,GAMMel93}). This set may vary with $y$, but in a number of important geometric situations in the edge (complete) setting it does not, see for example Mazzeo and Melrose \cite{GAMMaMe87}, Mazzeo and Phillips \cite{GAMMaPhil90}, Epstein, Melrose, and Mendoza \cite{GAMEpRBMMe91} to name but a few. The set
\begin{equation*}
\GAMspec_e(A)=\GAMset{(\sigma,y)\in \GAMC\times \GAMY:\sigma\in \GAMspec_{b,y}(A)}
\end{equation*}
is the edge spectrum of $A$. 

Thus, with the assumed ellipticity of $A$, \eqref{BeginExpansion} holds for sure in
\begin{equation*}
\GAMset{(\sigma,p)\in \GAMC\times \GAMN:\Im \sigma>\gamma-1,\ (\sigma,\wp(p))\notin \GAMspec_e(A)}.
\end{equation*}
Using this in the right-hand side of \eqref{BeginExpansion} gives now information about the meromorphic structure of $\widehat u(\sigma)$ in $\Im\sigma>\gamma-2$. Iterating, one gets information on $\Im\sigma>\gamma-m$. 

The caveat is that the right hand side has to be treated as a distribution (at least in the $y$ variable) so one has to proceed with much care.

Clearly, the meromorphic invertibility of $\GAMbPhat_y(\sigma)$ for $\sigma$ in the set 
\begin{equation*}
\Sigma=\GAMset{\sigma\in \GAMC:\gamma-m<\Im\sigma<\gamma}
\end{equation*}
for each $y\in\GAMY$ does play an important role. If $K$ is a Hilbert space and $V\subset \GAMC$ is open, we write $\GAMMero(V,K)$ for the space of meromorphic $K$-valued functions on $V$ and $\GAMHol(V,K)$ for the subspace of holomorphic elements. Thus $f\in \GAMMero(V,K)$ if there is, for each $\sigma_0\in V$, a number $\mu_0\in \GAMNN_0$ such that $\sigma\mapsto (\sigma-\sigma_0)^{\mu_0}f(\sigma)$ is holomorphic near $\sigma_0$. With this notation we have that the holomorphic (polynomial) family $\GAMbPhat_y$ defines a map
\begin{equation*}
\GAMbPhat_y:\GAMMero(\Sigma;H^m(\GAMZ_y;E_{\GAMZ_y}))\to \GAMMero(\Sigma;L^2(\GAMZ_y;F_{\GAMZ_y}))
\end{equation*}
which in turn gives a map
\begin{multline*}
[\GAMbPhat_y]:\GAMMero(\Sigma;H^m(\GAMZ_y;E_{\GAMZ_y}))/\GAMHol(\Sigma;H^m(\GAMZ_y;E_{\GAMZ_y}))\\\to \GAMMero(\Sigma;L^2(\GAMZ_y;F_{\GAMZ_y}))/\GAMHol(\Sigma;L^2(\GAMZ_y;F_{\GAMZ_y}))
\end{multline*}
It is the kernel of this last map that is of interest, as $y$ varies in $\GAMY$. This space is more conveniently expressed with the space of singular parts of its elements. The latter form a finite-dimensional space of meromorphic $H^m(\GAMZ_y;E_{\GAMZ_y})$-valued functions on $\GAMC$ with poles in $\GAMspec_{b,y}(A)\cap \Sigma$. To get an analytic hold on these elements, we note that the singular part of an element $[\widehat \tau]\in \ker [\GAMbPhat_y]$ represented by $\widehat \tau\in \GAMMero(\Sigma;H^m(\GAMZ_y;E_{\GAMZ_y}))$ can be obtained by computing the integral
\begin{equation*}
\ssp_\Omega(\widehat \tau)(\sigma)=\frac{\GAMim}{2\pi}\oint_{\partial\Omega}\frac{\widehat \tau(\zeta)}{\zeta-\sigma}\,d\zeta,\quad \ |\sigma|\gg 1
\end{equation*}
with the positive (counterclockwise) orientation, where $\Omega\Subset \Sigma$ is open, has smooth (or rectifiable) boundary and contains $\GAMspec_{b,y}(A)$. The integral is of course independent of the representative of $[\widehat \tau]$. Write
\begin{equation*}
\widehat \GAMtrb_y=\ssp(\ker[\GAMbPhat_y]).
\end{equation*}
If $\widehat \tau\in \widehat \GAMtrb_y$ and $\Gamma$ is a simple closed smooth (or rectifiable) curve surrounding $\GAMspec_{b,y}(A)\cap \Sigma$, then 
\begin{equation}\label{InversePseudoMellin}
\tau=-\frac{1}{2\pi}\oint_\Gamma x^{\GAMim\sigma}\widehat\tau(\sigma)\,d\sigma,
\end{equation}
with the positive orientation viewed as a section of $E_{\GAMZ_y^\wedge}$, is an element of $\ker \GAMA_y$ ($\tau$ is essentially the inverse Mellin transform of $\widehat \tau$).  Indeed,
\begin{align*}
x^m\GAMA_y\tau&=-\frac{1}{2\pi}\oint_\Gamma x^{\GAMim\sigma}\GAMbPhat(\sigma)u(\sigma)\,d\sigma\\
&=0
\end{align*}
since $\GAMbPhat(\sigma)u(\sigma)$ is entire. Evidently $\tau$ has the from
\begin{equation}\label{FiberOfTau}
\tau=\sum_{\substack{\sigma\in \Sigma\\ \sigma\in \GAMspec_{b,y}(A)}}\sum_{\ell=0}^{N_\sigma} \tau_{\sigma \ell}\, x^{\GAMim\sigma}\log^\ell x.
\end{equation}
The coefficients $\tau_{\sigma\ell}$ are sections (\emph{a fortiori} smooth) of $E_{\GAMZ_y}$. We let $\GAMtrb_y$ be the space of all elements obtained from $\widehat\GAMtrb_y$ as just described. 

\begin{theorem}[{\cite[Theorem {3.2}]{GAMKrMe12a}}]\label{TraceBundle}
Assume
\begin{equation}\label{FiniteSpecb}
\GAMset{(\sigma,y)\in \GAMC\times\GAMY:\Im\sigma=\gamma,\ \gamma-m}\cap \GAMspec_e(A)=\emptyset.
\end{equation}
Define
\begin{equation*}
\GAMtrb=\bigsqcup_{y\in \GAMY} \GAMtrb_y,\quad \pi:\GAMtrb\to \GAMY\text{ the canonical map}.
\end{equation*}
Then $\GAMtrb\to\GAMY$ is a smooth vector bundle. A smooth section of $\GAMtrb$ is a map
\begin{equation*}
\GAMY\ni y\mapsto \tau(y)\in \GAMtrb_y
\end{equation*}
which viewed as a section of $E$ over $\GAMopen \GAMN^\wedge$ is smooth in the usual sense.
\end{theorem}

We call $\GAMtrb$ the trace bundle of $A$. In this note we will always assume that \eqref{FiniteSpecb} holds.

\begin{example}
We pause to illustrate some of the ideas. Let $\GAMY$ be a closed orientable surface, $\GAMZ$ an arbitrary closed manifold, $L^+\to \GAMY$ a nontrivial complex line bundle, $L^-$ its dual, and
\begin{equation*}
\Phi=\begin{bmatrix}\Phi_{11}&\Phi_{12}\\\Phi_{21}&\Phi_{22}\end{bmatrix}
:\begin{matrix} L^+\\\oplus\phantom{{}^+}\\ L^-\end{matrix}\to \begin{matrix} L^+\\\oplus\phantom{{}^+}\\ L^-\end{matrix}
\end{equation*}
some vector bundle homomorphism covering the identity $\GAMY\to\GAMY$. The homomorphisms $\Phi_{11}$ and $\Phi_{22}$ are required to be isomorphisms. Let $\GAMM=\GAMlbra 0,1\GAMrpar\times \GAMY\times\GAMZ$ and let $E$ be the pullback of $L^+\oplus L^-$ to $\GAMM$, a trivial rank 2 bundle: $E=\pi^*L^+\oplus \pi^* L^-$, $\pi:\GAMM\to\GAMY$ the projection. Let $\Delta_{L^+}$ be a Laplacian on sections of $L^+$. For example, fix a Riemannian metric on $T^*\GAMY$, a connection and a Hermitian metric on $L^+$ and let $\Delta_{L^+}$ be the Bochner Laplacian, $\Delta_{L^+}=\nabla^*\nabla$. Let $\Delta_{L^-}$ be a Laplacian on sections of $L^-$. Exploiting the product structure of $\GAMM$ we let these Laplacians act on the factors of $E$ in the natural way. Next we pick a Laplacian $\Delta_\GAMZ$ (acting on functions on $\GAMZ$) and let $Q_\GAMZ=\Delta_\GAMZ+c$, $c\geq 0$, act on sections of $\pi^*L^+$ or $\pi^*L^-$ in the canonical way, again exploiting the product structure of $\GAMM$ and the pull-back nature of these bundles; we specify the constant $c$ later. Finally, writing $x$ for the standard coordinate in $\GAMlbra0,1\GAMrpar$ we define
\begin{equation*}
P=
\begin{bmatrix}(x D_x)^2&0\\0&(x D_x)^2\end{bmatrix}+
x^2\begin{bmatrix}\Delta_{L^+}&0\\0&\Delta_{L^-}\end{bmatrix}+
\begin{bmatrix}Q_\GAMZ\Phi_{11}&\Phi_{12}\\\Phi_{21}&Q_\GAMZ\Phi_{22}\end{bmatrix}
\end{equation*}
acting on $C^\infty(\GAMM;E)$. This is an elliptic edge operator with respect to the obvious boundary fibration of $\GAMM$, acting on sections of $\pi^*L^+\oplus\pi^*L^-$, so $A=x^{-2}P\in x^{-2}\GAMDiff^2_e(\GAMM;E)$ is an elliptic wedge operator. We compute its boundary spectrum. The indicial family of $P$ is
\begin{equation*}
\GAMbPhat(\sigma)=\begin{bmatrix}\sigma^2&0\\0&\sigma^2\end{bmatrix}+
\begin{bmatrix}Q_\GAMZ \Phi_{11}&\Phi_{12}\\\Phi_{21}&Q_\GAMZ\Phi_{22}\end{bmatrix}.
\end{equation*}
Let $\GAMset{\psi_k}$ be a complete orthonormal system of eigenvectors of $Q_\GAMZ$, $\GAMset{\lambda_k^2}$ the corresponding eigenvalues. The $\lambda_k$ are nonnegative and assumed to form a nondecreasing sequence. Let $\nu^+$ be a frame of $L^+$ in a neighborhood $U\subset \GAMY$ of some  $y_0\in\GAMY$, let $\nu^-$ be the dual frame.  Keeping the same notation for the lifted frames, they make up a frame for $E$ over $\pi^{-1}(U)$. In terms of the resulting frame, a section of $E_\GAMN$ over $\wp^{-1}(U)$ has the form
\begin{equation*}
\sum_{k=0}^\infty\psi_k(z)\begin{bmatrix} c_k^+(y)\\[3pt]c_k^-(y)\end{bmatrix}
\end{equation*}
Using formulas such as $\Phi_{12}(\nu^-)=\varphi_{12}\nu^+$ with smooth $\varphi_{ij}:U\to\GAMC$ we get 
\begin{equation}\label{bPhatAsMatrix}
\sum_{k=0}^\infty \psi_k(z)
\begin{bmatrix}\sigma^2 +\lambda_k^2\varphi_{11}(y)& \varphi_{12}(y)\\[3pt] \varphi_{21}(y)&\sigma^2+\lambda_k^2\varphi_{22}(y) \end{bmatrix}
\begin{bmatrix} c_k^+(y)\\[3pt]c_k^-(y)\end{bmatrix}
\end{equation}
for $\GAMbPhat(\sigma)u$ with respect to the same frame. The boundary spectrum at $y$ consists of the roots $\sigma$ of all polynomials
\begin{equation*}
\sigma^4+\lambda_k^2(\varphi_{11}(y)+\varphi_{22}(y))\sigma^2+\lambda_k^4\varphi_{11}(y)\varphi_{22}(y)-\varphi_{12}(y)\varphi_{21}(y),\quad k=0,1,\dots
\end{equation*}
The functions $\varphi_{11}$, $\varphi_{12}$ are globally defined and independent of choice of frames. The product $\varphi_{12}\varphi_{21}$ is independent of the choice of frame so it is also globally defined. Evidently, as $y$ varies, these roots can have very complicated behavior.

We illustrate the simplest possibilities assuming $\Phi_{21}=0$, in which case
\begin{equation*}
\GAMspec_{b,y}(A)=\GAMset{\pm\GAMim \lambda_k\sqrt{\varphi_{11}(y)},\ \pm\GAMim \lambda_k\sqrt{\varphi_{22}(y)}:k\in \GAMNN_0}.
\end{equation*}
Take $\gamma=1$. Pick $c>0$ in the definition of $Q_\GAMZ$ so $\lambda_0>0$. Suppose $1/\lambda_1^2<\sup_\GAMY|\varphi_{jj}|<1/\lambda_0^2$. Then \eqref{FiniteSpecb} holds and the part of the $b$-spectrum that matters, the part in $\Sigma=\GAMset{\sigma\in \GAMC:-1<\Im\sigma<1}$, consists of the points
\begin{equation*}
\pm\GAMim \lambda_0\sqrt{\varphi_{11}(y)},\ \pm\GAMim \lambda_0\sqrt{\varphi_{22}(y)}.
\end{equation*}
Each of these is locally well defined in a manner that gives a locally smooth function of $y$, since the $\varphi_{jj}$ vanishes nowhere. That the roots can be arranged to depend locally smoothly on $y$ cannot be guaranteed if $\varphi_{12}\varphi_{21}\ne 0$.
\begin{enumerate}
\item
If $\varphi_{11}(y)-\varphi_{22}(y)\ne0$ for all $y$ then these four roots are all different from each other for each $y$ , and $\widehat\GAMtrb_y$ is the span of
\begin{equation*}
\frac{\psi_1(z)}{\sigma\mp\GAMim\lambda_0\sqrt{\varphi_{11}(y)}} \nu^+,\qquad
\frac{\psi_1(y)}{\sigma\mp\GAMim\lambda_0\sqrt{\varphi_{22}(y)}}\Big(\dfrac{-\varphi_{12}}{\lambda_0^2(\varphi_{11}-\varphi_{22})}\nu^++\nu^-\Big).
\end{equation*}
The fibers of $\GAMtrb$ are $4$-dimensional, the trace bundle is isomorphic to $L^+\oplus L^+\oplus L^-\oplus L^-$ which is trivial since $\GAMchern_1(L^+)+\GAMchern_1(L^-)=0$, while neither individual summand, nor $L^+\oplus L^+$ and $L^-\oplus L^-$ is trivial since $2\GAMchern_1(L^+)\ne 0$. In particular, if $\GAMspec_{b,y}(A)$ is independent of $y$ (which happens when $\varphi_{11}$ and $\varphi_{22}$ are constant), the part of $\GAMtrb$ associated with a single pole in $\GAMspec_e(A)\cap \Sigma\times\GAMY$ is not trivial despite the fact that $E$ itself is a trivial bundle.
\item
Suppose now that $\varphi_{11}(y_0)=\varphi_{22}(y_0)$ for some $y_0$. For $y$ near $y_0$ define  $\sigma^\pm_1(y)=\pm \GAMim\lambda_0\sqrt{\phi_{11}(y)}$, likewise $\sigma^\pm_2(y)$. Then $\widehat\GAMtrb_{y_0}$ is spanned by 
\begin{equation*}
\chi^{\pm}_1(y_0)=\frac{1}{2\sigma^\pm_1(y_0)(\sigma-\sigma^\pm_1(y_0))}\begin{bmatrix}1\\0\end{bmatrix}
\end{equation*}
and
\begin{multline*}
\chi^\pm_2(y_0)=\bigg(\frac{1}{(\sigma-\sigma^\pm_1(y_0))^2}-\frac{1}{\sigma^\pm_1(y_0)(\sigma-\sigma^\pm_1(y_0))}\bigg)\!\!\begin{bmatrix}1\\0\end{bmatrix}\\ + \frac{{2\sigma^\pm_1(y_0)^2}}{\sigma^\pm_1(y_0)(\sigma-\sigma^\pm_1(y_0))}\!\!\begin{bmatrix}0\\1\end{bmatrix}
\end{multline*}
In the notation for the $\chi^\pm_j(y_0)$ we are taking advantage of the fact that $\varphi_{11}(y_0)=\varphi_{22}(y_0)$ and identify $\nu^+$ and $\nu^-$ with the respective column matrices. These formulas are obtained by applying the inverse of
\begin{equation*}
M(\sigma,y)=\begin{bmatrix}\sigma^2 +\lambda_0^2\varphi_{11}(y)& \varphi_{12}(y)\\[3pt] 0&\sigma^2+\lambda_0^2\varphi_{22}(y) \end{bmatrix}
\end{equation*}
(the matrix in \eqref{bPhatAsMatrix} with $k=0$ and $\varphi_{21}=0$) at $y=y_0$ to the column with components $f^+$, $f^-$ and then computing the singular part of the resulting expression at $\sigma=\sigma^\pm_1(y_0)$, which gives $f^+\chi^\pm_1(y_0)+f^-\chi^\pm_2(y_0)$. One now obtains a frame for $\widehat\GAMtrb$ near $y_0$ by taking the singular part of each of $M^{-1}(\sigma,y)M(\sigma,y_0)\chi^\pm_j$:
\begin{equation*}
\chi^\pm_j(\sigma,y) = \ssp_{\Omega}(\chi^\pm_j(\cdot, y_0))(\sigma)
\end{equation*}
where $\Omega$ is a disk containing the $\sigma^\pm_{jj}(y)$ in its interior ($y$ is kept in a neighborhood $V$ of $y_0$). Using the components of these four vectors as coefficients one rewrites the result in terms of $\nu^+$ and $\nu^-$. The resulting expressions for the $\chi^\pm_j(\sigma,y)$ are smooth in the complement of $(\Sigma\times V)\cap \GAMspec_e(A)$, and so applying \eqref{InversePseudoMellin} to each of them gives a smooth local section of $\GAMtrb$ over $V$, altogether making up a local frame.
\end{enumerate}
\end{example}

In general, the construction of local frames of $\GAMtrb$ near some $y_0$ as described at the end of the example yields pointwise bases for $\GAMtrb_y$ that are smooth in the sense described in Theorem~\ref{TraceBundle}. In \cite{GAMKrMe12a} we also prove that if $A^\star$ is the formal adjoint of $A$ and $\GAMtrb^\star$ its trace bundle, then, taking $\gamma=m/2$ for convenience, 
\begin{equation}\label{Pairing}
\GAMtrb_y\times \GAMtrb^\star_y \ni (u,v)\mapsto [u,v]^\flat_y = (\GAMA \omega u,\omega v)_{x^{-m/2}L^2_b}-(\omega u,\GAMA^\star\omega v)_{x^{-m/2}L^2_b}\in \GAMC
\end{equation}
is nondegenerate and gives a smooth Hermitian pairing of $\GAMtrb$ and $\GAMtrb^\star$. This is Theorem 5.3 of \cite{GAMKrMe12a}. The pairing is independent of the specific cut-off function $\omega$, but one needs to be included because the $L^2$ spaces are over $\GAMZ_y^\wedge$. What is behind nondegeneracy is the general fact that if $A$ is an arbitrary elliptic operator on some open $\GAMopen \GAMM$, then the Hilbert space adjoint of $A$ with its minimal domain is $A^\star$ with its maximal domain. The smoothness of the pairing follows from writing the pointwise pairing as a contour integral of Mellin transforms (as was in fact done in the just cited paper).

The proof in \cite{GAMKrMe12a} of Theorem~\ref{TraceBundle} above proceeds in two main steps. First we show the existence, in a neighborhood of each $y_0\in \GAMY$, of a system of sections $\tau_j$ that are smooth in the sense of the theorem and are a pointwise basis of each fiber. The existence of such local systems of solutions was also proved by Costabel and Dauge \cite{GAMCostabelDauge92} and Schmutzler \cite{GAMSchm92} using different methods. We then prove that two frames $\tau_j$, $\tau_j'$ are related by smooth transition functions by observing that if $\tau_j^\star$ is a local frame for $\GAMtrb^*$, then, first 
\begin{equation*}
\tau_j'=\sum_k a_{kj}\tau_k
\end{equation*}
for some functions $a_{kj}$, trivially since the $\tau_j$ give bases pointwise, and then that the $a_{kj}$ must be smooth because they satisfy the system of equations
\begin{equation*}
[\tau_j',\tau_\ell^\star]^\flat_y=a_{kj}(y)[\tau_k,\tau_\ell^\star]^\flat_y
\end{equation*}
in which the matrices with components $[\tau_j',\tau_\ell^\star]^\flat_y$ and $[\tau_k,\tau_\ell^\star]^\flat_y$ are smooth, since the $\tau_j'$, $\tau_k$, and $\tau_\ell^\star$ are smooth, and the second matrix is invertible by the nondegeneracy of the pairing. It follows that the set of frames which are smooth in the sense of Theorem~\ref{TraceBundle} admits smooth transition functions, so the condition defines a smooth structure for the trace bundle in the usual sense. 

We close this section with one last observation whose relevancy will become apparent in each of the next two sections. The fiber of $\GAMtrb$ at $y$ consists of elements of the form \eqref{FiberOfTau} in the kernel of $\GAMA_y$. Since $\GAMA_y$ and $\nabla_{x\partial_x}$ commute, the latter defines a bundle homomorphism which we shall denote by
\begin{equation*}
x\partial_x:\GAMtrb\to\GAMtrb.
\end{equation*}
This homomorphism is smooth because at the level of $\widehat\GAMtrb_y$ it is just multiplication by $\GAMim\sigma$. The eigenvalues on the fiber $\GAMtrb_y$, the numbers $\GAMim\sigma$ with $\sigma\in \GAMspec_{b,y}(A)\cap \Sigma$, generally vary with the fiber, as will the Jordan canonical form of $x\partial_x$.

%%%%%%%%%%%%%%%%%%%%%%%%%%%%%%%%%%%%%%%%%%%%%%%%%
\section{Elliptic systems of variable order}\label{sec-variableorder}

To motivate the results described in this section, mostly coming from \cite{GAMKrMe12b}, it is useful to follow the construction of the trace bundle in the case of a regular elliptic differential operator on a manifold with boundary. 

\begin{example}\label{ClassicalTrace}
Suppose $A$ is such an operator on $\GAMM$. In local coordinates $x,y_j$ near a point of $\GAMN=\GAMY$ (each point of $\GAMN$ is a fiber of the boundary fibration, so there are no $z_\mu$),
\begin{equation*}
A=\sum_{k+|\alpha|\leq m}a_{k\alpha}(x,y) D_x^kD_y^\alpha
\end{equation*}
with smooth $a_{k\alpha}$ up to $x=0$. Using 
\begin{equation*}
x^m D_x^k D_y^\alpha=x^{m-k-|\alpha|}p_k(xD_x+\GAMim|\alpha|)(xD_y)^\alpha
\end{equation*}
one gets
\begin{equation*}
x^mA=\sum_{k+|\alpha|\leq m}a_{k\alpha}(x,y)x^{m-k-|\alpha|} p_k(xD_x+\GAMim|\alpha|)(xD_y)^\alpha,
\end{equation*}
so
\begin{equation*}
\GAMbPhat_y(\sigma)=a_{m0}(0,y)p_m(\sigma)
\end{equation*}
which translates to $\GAMbP_y=a_{m0}(0,y)p_m(xD_x)$, therefore
\begin{equation*}
\GAMA_y=a_{m0}(0,y)D_x^m.
\end{equation*}
Taking $\gamma=1/2$, the relevant strip $\Sigma$ in the complex plane is $1/2-m<\Im \sigma<1/2$, the set of poles of $\GAMbPhat_y(\sigma)^{-1}$ is
\begin{equation*}
\GAMspec_{b,y}(A)=\GAMset{0,-\GAMim,-2\GAMim,\dots,-(m-1)\GAMim},
\end{equation*}
all poles are simple, and the elements of $\widehat\GAMtrb_y$ have the form
\begin{equation*}
\widehat \tau=\sum_{\ell=0}^{m-1} \frac{\tau_\ell}{\sigma+\GAMim \ell}
\end{equation*}
where $\tau_\ell\in E_{y}$. Thus $\GAMtrb_y$ consists of all polynomials
\begin{equation*}
\tau=-\frac{1}{2\pi}\sum_{\ell=0}^{m-1} \oint_\Gamma x^{\GAMim\sigma}\frac{\tau_\ell}{\sigma+\GAMim \ell}\,d\sigma=-\frac{1}{2\pi}\sum_{\ell=0}^{m-1}\tau_\ell x^\ell
\end{equation*}
as functions on $\GAMR_+$ with values in the fiber $E_y$ of $E$. These are exactly the elements in the kernel of $\GAMA_y$ since $a_{m0}(0,y)$ is invertible by ellipticity. Therefore, by our definition, the trace bundle of $A$ is
\begin{equation*}
\GAMtrb=\bigsqcup_{y\in \GAMY} \Big\{\sum_{\ell=0}^{m-1}\tau_\ell x^\ell: \tau_\ell\in E_y\Big\},
\end{equation*}
the direct sum of $m$ copies of $E_\GAMY$. It is a particularity of regular elliptic operators of order $m$ on sections of a bundle $E$ that they all share the same trace bundle. As expected, the operator $x\partial_x$ discussed in the last paragraph of the previous section acts on $\GAMtrb$. Its eigenvalues are the numbers $0,1,\dots,m-1$, the eigenvectors in the fiber over $y\in \GAMY$ corresponding to the eigenvalue $\ell$ are the monomials $\tau_\ell x^\ell$, $\tau_\ell\in E_y$. The trace bundle splits globally into the direct sum of $m$ subbundles, each isomorphic to $E$.

\end{example}

Connecting the above example with a boundary problem
\begin{equation*}
\begin{cases}
Au=f&\text{ in }\GAMopen\GAMM\\
B\gamma u=g&\text{ on }\GAMY
\end{cases}
\end{equation*}
for $A$ (where the solution $u$ is sought in $H^m(\GAMM;E)$), note that the classical traces $\gamma_s(u)=D_x^s u\big|_\GAMY$, assembled into the Taylor polynomial
\begin{equation*}
\gamma_A(u)=\sum_{\ell=0}^{m-1}\frac{\GAMim^\ell}{\ell!}\gamma_\ell(u)\, x^\ell.
\end{equation*}
of $u$, yield a section of the trace bundle of $A$. The relation between the regularity of the coefficients $\gamma_\ell(u)$ and the power $x^\ell$, namely that the component of $\gamma_A(u)$ in the eigenspace of $x\partial_x$ with eigenvalue $\ell$ lies in $H^{m-\ell-1/2}(\GAMY;\GAMtrb)$, is not to be viewed as an accident but as an expression of a tight link between these concepts. This assertion will be fully justified by the results to be described in the next section. 

In the case of a general elliptic wedge operator $A$, the fiberwise action of $x\partial_x$ on its trace bundle will generically have eigenvalues and Jordan canonical form varying with the base point. This makes the issue of regularity of boundary values rather more complicated. We deal with this in \cite{GAMKrMe12b} as a problem independently of the motivating example. In the rest of the section we present some of the ideas going into that paper, mostly from the local point of view, rarely yielding to expressing things globally.

\medskip
On account that the results are independent of those of \cite{GAMKrMe12a} which were described in the previous section, we view $\pi:\GAMtrb\to \GAMY$ as some vector bundle and consider an arbitrary smooth endomorphism $\GAMagen:\GAMtrb\to \GAMtrb$ on which no conditions are placed. We continue to assume that $\GAMY$ is compact although this is not necessary in the general theory. We fix some Hermitian metric on $\GAMtrb$, not necessarily related to $\GAMagen$, and a smooth density $\GAMm_\GAMY$. With these we define $L^2(\GAMY;\GAMtrb)$.

\medskip
The following example motivates the next step:

\begin{example}
Let $a: \GAMC^m\to \GAMC^m$ be diagonal with entries $\ell+1/2$, $\ell=0,\dotsc,m-1$. For $\eta\in \GAMR^q$ define $\GAMjb{\eta}=(1+|\eta|^2)^{1/2}$ and set $\GAMjb{\eta}^{m-a}=\exp\big(\log\GAMjb{\eta}(m-a)\big)$. Then $\GAMjb{\eta}^{m-a}$ is diagonal with entries $\GAMjb{\eta}^{m-\ell-1/2}$ and if $u\in C_c^{-\infty}(\GAMR^q;\GAMC^m)$ is such that
\begin{equation*}
\frac{1}{(2\pi)^q}\int_{\GAMR^q} e^{\GAMim y\cdot \eta}\GAMjb{\eta}^{m-a}\,\widehat u(\eta)\,d\eta
\end{equation*}
belongs to $L^2(\GAMR^d;\GAMC^m)$, then the components $u_0,\dotsc,u_{m-1}$ satisfy 
\begin{equation*}
u_\ell\in H^{m-\ell-1/2}(\GAMR^q).
\end{equation*}
\end{example}

Returning to the general case, define $\varrho^{\GAMagen}$ for $\varrho>0$ as expected:
\begin{equation*}
\varrho^{\GAMagen(y)}=\frac{\GAMim}{2\pi}\oint_{\Gamma_y} \varrho^\sigma(\GAMagen(y)-\sigma)^{-1}\,d\sigma,\quad \varrho>0,
\end{equation*}
where $\Gamma_y$ is a positively oriented contour enclosing $\GAMspec(\GAMagen(y))$. This defines a smooth isomorphism $\GAMtrb\to\GAMtrb$. Finally, let $g$ be a Riemannian metric on $\GAMY$, define $\GAMjb{\pmb\eta}=(1+g(\pmb\eta,\pmb\eta))^{1/2}$, then $p(\pmb\eta)=\GAMjb{\pmb\eta}^{\GAMagen(y)}$, for $\pmb\eta\in T_y^*\GAMY$. We will prove in a moment that $p$ is a symbol in the class $S^\infty_{1,\delta}$ for $\delta>0$ arbitrarily small, locally in any sufficiently small neighborhood of any point of $\GAMY$. The order is locally bounded (by our assumption of compactness of $\GAMY$, also globally bounded).

To define a pseudodifferential operator with $p$ as principal symbol, construct a global parametrization of the conormal to the diagonal $\Delta_\GAMY\subset \GAMY\times\GAMY$ as in \cite[pg. 381]{GAMGuiSt77}. Namely, view $T\GAMY$ as the normal bundle to $\Delta_\GAMY$: $N\Delta_\GAMY=\GAMset{(v,-v):v\in T\GAMY}$, let $V$ be a neighborhood of the zero section of $T\GAMY$ on which $\exp:V\to \GAMY\times\GAMY$ is a diffeomorphism onto its image $W$. If $(y,y')\in W$, let $(v,-v)$, $v\in T_{y''}\GAMY$ be the unique element mapped to $(y,y')$ by $\exp$ and define $\varphi(y,y',\pmb\eta)=\langle\pmb\eta,v\rangle$ for $\pmb\eta\in T^*_{y''}\GAMM$. Finally, pick a smooth homogeneous fiberwise density $\nu$ on $T^*\GAMY$ (in local coordinates $(y,\eta)$, $d\nu=f(y)|d\eta_1\wedge\dots\wedge d\eta_q|$, $f$ smooth and positive). With these ingredients define $\Lambda^{\GAMagen}$ through its Schwartz kernel:
\begin{equation*}
K_{\Lambda^{\GAMagen}}(y,y')=\frac{1}{(2\pi)^q}\int e^{\GAMim\varphi(y,y',\pmb\eta)}\sum_\alpha \chi_\alpha(y,y') p_\alpha(y,y',\eta)\,d\nu
\end{equation*}
where the $\chi_\alpha$ are carefully chosen smooth functions with $\sum\chi_\alpha$ supported in $W$, equal to $1$ in a neighborhood of $\Delta_\GAMY$, and the $p_\alpha$ are constructed using $\GAMjb{\eta}_y^\GAMagen$ using certain adapted trivializations as we explain below. We then define 
\begin{equation*}
H^{\GAMagen+s}(\GAMY;E)=\GAMset{u\in C^{-\infty}(\GAMY;\GAMtrb): \Lambda^{\GAMagen+s} u\in L^2(\GAMY;E)}
\end{equation*}
for real $s$ using the endomorphism $\GAMagen+s\GAMId$ of $\GAMtrb$.

\medskip
It remains to explain the choice of functions $\chi_\alpha$ and what the $p_\alpha$ are. The purpose is to gain some control on the behavior of $\GAMagen$ through careful localization as follows. Fix $\delta\in (0,1)$. Given $y_0\in \GAMY$ choose, for each eigenvalue $\sigma_\ell$ of $\GAMagen(y_0):\GAMtrb_{y_0}\to \GAMtrb_{y_0}$ a number $\delta_\ell \in (0,\delta)$ small enough that the closures of the disks $D_\ell=D(\sigma_\ell,\delta_\ell/2)\subset \GAMC$ are disjoint. Then there is a neighborhood $U$ of $y_0$ such that $\GAMspec(\GAMagen(y))\subset \bigcup D_\ell$ if $y\in U$. By way of the projections
\begin{equation*}
\Pi_{\ell,y}=\frac{\GAMim }{2\pi}\oint_{\partial D_\ell}(\GAMagen(y)-\sigma)^{-1}\,d\sigma,\quad y\in U,
\end{equation*}
one gets a decomposition
\begin{equation*}
\GAMtrb_U=\bigoplus_\ell \GAMtrb^\ell_U,\qquad \GAMtrb^\ell_U=\Pi_\ell\GAMtrb_U
\end{equation*}
of $\GAMtrb$ over $U$ into smooth subbundles, each invariant under $\GAMagen$; the eigenvalues of $\GAMagen|_{\GAMtrb^\ell}:\GAMtrb^\ell_U\to \GAMtrb^\ell_U$ cluster within $D(\sigma_\ell,\delta_\ell/2)$. Picking $U$ small enough allows us to assume additionally that the bundles $\GAMtrb^\ell_U$ (hence also $\GAMtrb_U$) are trivial. We refer to the above as a $\delta$-admissible decomposition of $\GAMtrb$ over $U$ (\cite[Definition 2.2]{GAMKrMe12b}). We trivialize $\GAMtrb_U$ though the trivializations of the $\GAMtrb^\ell_U$: Let $\phi:\pi^{-1}(U)\to U\times \GAMC^r$ be such an adapted  trivialization ($r=\GAMrank \GAMtrb$), define $a_\phi:U\times\GAMC^r\to U\times\GAMC^r$ by $a_\phi=\phi\circ \GAMagen\circ \phi^{-1}$. Finally, assume that $U$ is also the domain of a local chart of $\GAMY$. 

We now pick an open cover $\GAMset{U_\alpha}$ of $\GAMY$ consisting of open sets as just described and build up the $\chi_{\alpha}\in C_c^\infty(U_\alpha\times U_\alpha)$ from a partition of unity near the diagonal subordinate to the cover $\GAMset{U_\alpha\times U_\alpha}$ of $\Delta_\GAMY$ so that $\sum \chi_\alpha=1$ near $\Delta_\GAMY$.  Next, with adapted trivializations $\phi_\alpha$ we let
\begin{equation*}
p_\alpha(y,y',\eta)=\phi_\alpha(y)^{-1}\GAMjb{\eta}_{y}^{a_{\phi_\alpha}(y)}\phi_\alpha(y'),
\end{equation*}
so that $p_\alpha(y,y',\eta):\GAMtrb_{y'}\to\GAMtrb_y$. Implicit in this definition is that we are using parallel transport on $\GAMtrb|_{U_\alpha}$ with respect to a (flat) connection adapted to a $\delta$-admissible decomposition. See Safarov \cite{GAMSaf97} or Pflaum \cite{GAMPfl98} for a systematic analysis of the role of connections in the definition of standard pseudodifferential operators acting on vector bundles. 

\medskip
It remains to show that the $p_\alpha$ are symbols. In the following lemma we let $\GAMjb{\eta}_y^2=1+\sum g^{ij}(y)\eta_i\eta_j$ with smooth, positive definite $g^{ij}$. We will drop the reference to $\phi$ from the notation. 

\begin{lemma}[{\cite[Lemma 3.4]{GAMKrMe12b}}]\label{StandardAction} The function $(y,\eta)\mapsto \GAMjb{\eta}_y^{a(y)}$ is, in each open subset $V\Subset U$, a symbol in the H\"ormander class $S^M_{1,\delta}$ for some $M$ depending on $V$.
\end{lemma}

The proof of the lemma consists of establishing that
\begin{equation*}
D_y^\alpha\partial_\eta^\beta \GAMjb{\eta}_y^{a(y)} = p_{\alpha \beta}(y,\eta)\GAMjb{\eta}_y^{a(y)}
\end{equation*}
with $\|p_{\alpha \beta}(y,\eta)\|\leq C\GAMjb{\eta}_y^{-|\beta|+\delta}$ (one can take $\delta=0$ in the estimate if $\alpha=0$), then observing that
\begin{equation*}
\|\GAMjb{\eta}_y^{a(y)}\|\leq C(1+|\eta|)^M,\quad (y,\eta)\in V\times\GAMR^q
\end{equation*}
for some $C$ if $M>\sup\GAMset{\Re\sigma:\sigma\in \GAMspec(a(y)),\ y\in V}$. 

To get the estimate for $p_{\alpha \beta}$, note first that since $\varrho\partial_\varrho \varrho^{a(y)}=a(y)\varrho^{a(y)}$, 
\begin{equation*}
\partial_{\eta_j}\GAMjb{\eta}_y^{a(y)}=\frac{\partial_{\eta_j}\GAMjb{\eta}_y}{\GAMjb{\eta}_y}a(y)\GAMjb{\eta}_y^{a(y)},
\end{equation*}
and by induction,
\begin{equation*}
\partial_\eta^\beta\GAMjb{\eta}_y^{a(y)}=p_\beta(y,\eta)\GAMjb{\eta}_y^{a(y)},
\end{equation*}
where $p_\beta(y,\eta)$ is a classical symbol of order $-|\beta|$. So we only need to deal with derivatives in $y$. Differentiating
\begin{equation*}
\GAMjb{\eta}_y^{a(y)} = \sum_j \frac{\GAMim}{2\pi}\oint_{\partial D_\ell}\GAMjb{\eta}_y^\sigma (a(y)-\sigma)^{-1}\,d\sigma
\end{equation*}
one gets that
\begin{equation*}
D_y^\alpha\GAMjb{\eta}_y^{a(y)} = \sum_\ell \frac{\GAMim}{2\pi}\oint_{\partial D_\ell}\GAMjb{\eta}_y^\sigma R_{\alpha}(y,\sigma)\,d\sigma
\end{equation*}
where $R_{\alpha}(y,\sigma)=\sum_{k=0}^{|\alpha|}\sigma^k R_{\alpha k}$ where $R_{\alpha k}$ is a sum of terms each being the product of a classical symbol of order zero followed by a product of at most $|\alpha|-k+1$ factors $(a(y)-\sigma)^{-1}$, each of these factors separated by a factor $\partial_y^{\alpha'}a(y)$. Since $\GAMdist(\GAMspec(a(y))\cap D_\ell,\partial D_\ell)$ is uniformly bounded from below when $y\in V$ (since $V\Subset U$), the norm of $R_\alpha(y,\sigma)$ is also uniformly bounded  when $y\in V$ and $\sigma\in \bigcup_\ell \partial D_\ell$. Now use 
\begin{equation*}
(D_y^\alpha\GAMjb{\eta}_y^{a(y)})\GAMjb{\eta}_y^{-a(y)} = \sum_\ell \Big(\frac{\GAMim}{2\pi}\Big)^2\int_{\partial D_\ell\times\partial D_\ell}\GAMjb{\eta}_y^{\sigma-\sigma'} R_\alpha(y,\sigma)(a(y)-\sigma')^{-1}\,d\sigma\,d\sigma'
\end{equation*}
to get the estimate 
\begin{equation*}
\|(D_y^\alpha\GAMjb{\eta}_y^{a(y)})\GAMjb{\eta}_y^{-a(y)}\|\leq C\GAMjb{\eta}^\delta
\end{equation*}
using that $\|\GAMjb{\eta}_y^{\sigma-\sigma'}\|\leq C(1+|\eta|)^{\delta}$ when $\sigma,\sigma'\in \partial D_\ell$. This completes the proof of the lemma. 

\medskip
In addition to defining the spaces $H^{s+\GAMagen}(\GAMY;\GAMtrb)$ to deal with the expected issue of varying regularity of traces of elements of the maximal domain of an elliptic wedge operator, one also needs to develop a theory of pseudodifferential operators able to deal with the varying regularity in order to express general boundary conditions. This of course is reminiscent of the relation between the Douglis-Nirenberg calculus \cite{GAMDouglisNirenberg}, see also Chazarain and Piriou \cite{GAMChazarainPiriou}, and boundary conditions in the classical situation. 

The ingredients going into the definition of such operators are a pair of vector bundles $\GAMtrb,\GAMtrt\to\GAMY$ endowed with endomorphisms $\GAMagen$, $\GAMbgen$. Symbol classes are defined locally on $\delta$-admissible domains common to both $\GAMagen$ and $\GAMbgen$ which are also domains of local charts. Passing to trivializations and local coordinates, we define (\cite[Definition 3.1]{GAMKrMe12b})
\begin{equation*}
S_{1,\delta}^{\mu}(U\times\GAMR^q;(\GAMC^{r},a),(\GAMC^{r'},b))
\end{equation*}
for any real $\mu$ to be the space of all $p(y,\eta) \in C^{\infty}(U\times\GAMR^q,\GAMHom(\GAMC^{r},\GAMC^{r'}))$ such that for every subset $K \Subset U$ and all $\alpha,\beta \in \GAMNN_0^q$ there exists a constant $C_{K,\alpha,\beta} > 0$ such that
\begin{equation*}
\|\langle \eta \rangle^{b(y)}\bigl(D_y^{\alpha}\partial_{\eta}^{\beta}p(y,\eta)\bigr)\langle \eta \rangle^{-a(y)}\| \leq C_{K,\alpha,\beta} \langle \eta \rangle^{\mu - |\beta| + \delta|\alpha|}
\end{equation*}
for all $(y,\eta) \in K\times\GAMR^q$. For example, the symbol in Lemma~\ref{StandardAction} belongs to the class $S_{1,\delta}^0(U\times\GAMR^q;(\GAMC^{r},a),(\GAMC^{r},0))$. It also belongs to $S_{1,\delta}^0(U\times\GAMR^q;(\GAMC^{r},0),(\GAMC^{r},-a))$, as can be seen by composing $D_y^\alpha\partial_\eta^\beta\GAMjb{\eta}_y^{a(y)}$ with $\GAMjb{\eta}^{-a}$ on the left rather than on the right.

The classes $S_{1,\delta}^{\mu}(U\times\GAMR^q;(\GAMC^{r},a),(\GAMC^{r'},b))$ have all the usual properties such as invariance under changes of coordinates, and asymptotic summability. Ellipticity of $p\in S_{1,\delta}^{\mu}(U\times\GAMR^q;(\GAMC^{r},a),(\GAMC^{r'},b))$ is defined as usual: the existence of $q\in S_{1,\delta}^{-\mu}(U\times\GAMR^q;(\GAMC^{r'},b),(\GAMC^{r},a))$ such that $pq=1$ modulo lower order terms.

Elements in the class just defined are locally in $S^M_{1,\delta}$ for large enough $M$, so one can use them to define pseudodifferential operators; write $\Psi^\mu_{1,\delta}(U;(\GAMC^{r},a),(\GAMC^{r'},b))$ for the space with symbols as just defined. The expected properties hold: invariance under changes of local coordinates and under compatible composition, the latter meaning
\begin{equation*}
\Psi^{\mu'}_{1,\delta}(U;(\GAMC^{r'},b),(\GAMC^{r''},c))\circ\Psi^{\mu}_{1,\delta}(U;(\GAMC^{r},a),(\GAMC^{r'},b))\subset \Psi^{\mu+\mu'}_{1,\delta}(U;(\GAMC^{r},a),(\GAMC^{r''},c)).
\end{equation*}
The reason for the good behavior at this level is that because the symbols are in a H\"ormander class, all the properties of these classes are inherited. In particular, one gets asymptotic expansions for the symbol of a composition, from which one can directly determine the validity of the assertion about compositions in our class.

Consequently the rough theory, including global definitions of classes of operators of varying order, follows largely as expected (even though some proofs had to be reworked completely) leading to the definition of the classes 
\begin{equation*}
\Psi^\mu_{1,\delta}(\GAMY;(\GAMtrb,\GAMagen),(\GAMtrt,\GAMbgen))
\end{equation*}
associated with pairs of vector-bundles-with-endomorphism. We single out a subclass having twisted homogeneous principal symbols (of order $\mu$), meaning that they are locally constructed using symbols with the property
\begin{equation*}
p(y,\varrho \eta) = \varrho^{\mu} \varrho^{-b(y)} p(y,\eta) \varrho^{a(y)}\quad \text{ for all $\varrho > 1$ and $|\eta|>1$}
\end{equation*}
modulo arbitrary symbols in $S^{\mu-1+\delta}_{1,\delta}(U;(\GAMC^r,a),(\GAMC^{r'},b))$. An operator $P$ in the class $\Psi^\mu_{1,\delta}(\GAMY;(\GAMtrb,\GAMagen),(\GAMtrt,\GAMbgen))$ has, as expected, a globally defined principal symbol
\begin{equation*}
\GAMsym(P):T^*\GAMY\GAMminus 0\to \GAMHom(\pi_\GAMY^*\GAMtrb,\pi_\GAMY^*\GAMtrt)
\end{equation*}
satisfying
\begin{equation*}
\GAMsym(P)(\varrho \pmb\eta) = \varrho^{\mu} \varrho^{-\GAMbgen(y)} \GAMsym(P)(\pmb\eta) \varrho^{\GAMagen(y)},\quad \varrho>0,\ \pmb\eta\in T^*\GAMY\GAMminus 0.
\end{equation*}

%%%%%%%%%%%%%%%%%%%%%%%%%%%%%%%%%%%%%%%%%%%%%%%%%
\section{Boundary value problems for first order elliptic wedge operators}%\label{sec-bvp}

In this section we describe some of the results contained \cite{GAMKrMe13}. As in that paper, we limit our discussion here to operators of order $1$, which allows us to circumvent a number of technical complications proper to the higher order case; these are being addressed in a forthcoming paper \cite{GAMKrMe15}. 

Henceforth $A$ is an elliptic element of $x^{-1}\GAMDiff^1_e(\GAMM;E,F)$. We will base our discussion on the spaces $x^{-1/2}L^2_b$, \GAMie, we set $\gamma=1/2$ (see the last paragraph of Sections~\ref{sec-setup} and~\ref{sec-considerations}). We aim at describing our specific approach to setting up boundary value problems for $A$, then establishing sufficient conditions for well-posedness of these problems.

Recall that the trace bundle of $A$ generally does depend on $A$ and that it carries a natural endomorphism $x\partial_x$. In brief, we need to describe a restriction, or trace, operator
\begin{equation*}
\gamma_A:\GAMDom_{\max}(A)\to C^{-\infty}(\GAMY;\GAMtrb),
\end{equation*}
to state the problem, then also ancillary objects to state sufficient conditions for well-posedness.

Our first observation is that if $u$ is a smooth section of $\GAMtrb$ and $\omega$ is a cut-off function, then $\omega u\in \GAMDom_{\max}(A)$ if $u\in C^\infty(\GAMY;\GAMtrb)$; this is Lemma 5.6 of \cite{GAMKrMe13}. Strictly speaking, $u$ is defined on $\GAMN^\wedge$ and is a smooth section of $E^\wedge$ over the interior of $\GAMN^\wedge$. Multiplication by $\omega$ produces an element of compact support which is transferred to $\GAMM$ via a tubular neighborhood map (and parallel transport). Furthermore, see Lemma~5.6, \GAMopcit, the map
\begin{equation*}
C^\infty(\GAMY;\GAMtrb)\ni u\mapsto \preP u=\omega u\in \GAMDom_{\max}(A)
\end{equation*}
is continuous.

We now use $\preP$ and the pairing \eqref{Pairing} to define $\gamma_{A^*}$. For details see Section 5, {\GAMopcit} Let $\preP^\star:C^\infty(\GAMY;\GAMtrb^\star)\to \GAMDom_{\max}(A^\star)$ be the corresponding map for the formal adjoint of $A$. Pick $v\in \GAMDom_{\max}(A^\star)$. Then
\begin{equation*}
C^\infty(\GAMY;\GAMtrb)\ni u\mapsto \langle\lambda_v,u\rangle=[\preP u,v]_A=(A\preP u,v)_{x^{-1/2}L^2}-(\preP u,A^\star v)_{x^{-1/2}L^2}\in \GAMC
\end{equation*}
is continuous, so it defines a distribution with values in the bundle $\GAMtrb^*\otimes|\GAMWedge|\GAMY$. Here $\GAMtrb^*$ is the dual bundle of $\GAMtrb$ with the opposite complex structure, not the trace bundle of $A^\star$, and $|\GAMWedge|\GAMY$ is the density bundle of $\GAMY$, which we trivialize using the density $\GAMm_\GAMY$. We get an antilinear map $v\mapsto \lambda_v$ which is converted into a linear map $v\to\gamma_{A^\star} v$ using the nondegenerate pairing $[\cdot,\cdot]^\flat$ to identify $\GAMtrb^*$ with $\GAMtrb^\star$, so that 
\begin{equation*}
[\preP u,v]_A=\int_\GAMY[u, \gamma_{A^\star}v]_y^\flat\,d\GAMm_\GAMY.
\end{equation*}
We now have:

\begin{theorem}[Theorem 5.11, \GAMopcit] The trace map is a continuous operator
\begin{equation*}
\gamma_A : \GAMDom_{\max}(A) \to H^{-x\partial_x -1/2}(\GAMY;\GAMtrb).
\end{equation*}
\end{theorem}

Note that if $A$ is a regular first order elliptic operator, then $x\partial x$ acts as $0$, since the traces are zeroth order polynomials by Example~\ref{ClassicalTrace}, so in this case the space $H^{-x\partial_x -1/2}(\GAMY;\GAMtrb)$ is the classical space of sections of $E_\GAMY$ of Sobolev regularity $-1/2$.

\medskip
With $\preP$ and $\gamma_A$ we can now define a replacement for the classical Sobolev space $H^1(\GAMM;E)$ in which solutions are sought when $A$ is a regular first order elliptic operator, see Section 8, {\GAMopcit} Namely, we let $H^1_\GAMtrb(\GAMM;E)$ be the completion of $C^\infty_\GAMtrb(\GAMM;E)=\dot C^\infty(\GAMM;E)+\preP C^\infty(\GAMY;\GAMtrb)$, the latter being a replacement of $C^\infty(\GAMM;E)$, with respect to the norm defined by
\begin{equation*}
\|u\|_{H^1_\GAMtrb}^2=\|u\|^2_A+\|\gamma_A u\|^2_{H^{1/2-x\partial_x}}.
\end{equation*}
We recall that $\|u\|^2_A=\|u\|^2+\|A u\|^2$ is the norm on $\GAMDom_{\max}(A)$ but that this norm is typically too weak to give a compact  embedding $\GAMDom_{\max}(A)\GAMembed x^{-1/2}L^2_b$, as illustrated in Example~\ref{BasicExample}.

\medskip
The following example shows how this trick works.

\begin{example}
Consider the situation of Example~\ref{BasicExample}, in which $\GAMM$ is the closed unit disk in $\GAMR^2$ and $\Delta$ the standard Laplacian. We claim that the norm defined by
\begin{equation*}
\|u\|_{H^2_\GAMtrb}^2 = \|u\|_\Delta^2+\|\gamma_0 u\|_{H^{3/2}(\partial \GAMM)}^2+\|\gamma_1 u\|_{H^{1/2}(\partial \GAMM)}^2
\end{equation*}
on $C^\infty(\GAMM)$, in which $\gamma_\ell=\partial^\ell u\big|_{\partial\GAMM}$, is equivalent to the $H^2$ norm on $\GAMM$. Indeed, on the one hand, the continuity of
\begin{equation*}
H^{2}(\GAMM)\xrightarrow{\gamma_\ell} H^{2-\ell-1/2}(\partial \GAMM),\quad \ell=0,1,
\end{equation*}
gives
\begin{equation*}
\|u\|_{H^2_\GAMtrb}^2\leq C\|u\|_{H^2(\GAMM)}^2
\end{equation*}
for some $C$. On the other hand,
\begin{equation*}
\int_\GAMM\|\nabla u\|^2\,d\lambda = - \int_\GAMM u\,\Delta u\,d\lambda - \int_{\partial\GAMM} u\, \frac{\partial\overline u}{\partial\nu}\, ds
\end{equation*}
and the Cauchy-Schwarz inequality give the reverse estimate; $\nu$ is the inward pointing normal, $d\lambda$ is Lebesgue measure and $ds$ arc-length measure. Thus, while the $\Delta$-graph norm \eqref{DeltaNorm} is not sufficiently strong to give compactness of the embedding $\GAMDom_{\max}(\Delta)\to L^2(\GAMM)$, the norm $\|u\|_{H^2_\GAMtrb}^2$ does. One gets automatically a split exact sequence
\begin{equation*}
0 \to H^2_0(\GAMM) \to H_\GAMtrb^2(\GAMM) \to H^{2-(x\partial_x+1/2)}(\partial\GAMM,\GAMtrb) \to 0
\end{equation*}
as we discussed in Section \ref{sec-considerations}.
\end{example}

In the general case we are discussing, the space $H^1_\GAMtrb(\GAMM;E)$ satisfies the basic requirements mentioned in Section~\ref{sec-considerations}. First, \eqref{SplitExactSequence} takes on the form
\begin{equation*}
0 \to \GAMDom_{\min}(A) \to H_\GAMtrb^{1}(\GAMM;E) \xrightarrow{\gamma_A} H^{1-(x\partial_x+1/2)}(\partial\GAMM,\GAMtrb) \to 0.
\end{equation*}
This is an exact sequence which splits though a continuous extension operator
\begin{equation*}
\GAMExt:H^{1-(x\partial_x+1/2)}(\partial\GAMM,\GAMtrb)\to H^1_\GAMtrb(\GAMM;E),
\end{equation*}
a left inverse of $\gamma_A$.  Second, the inclusion
\begin{equation*}
H^1_\GAMtrb(\GAMM;E)\GAMembed x^{-1/2}L^2_b(\GAMM;E)
\end{equation*}
is compact. We shall omit here the discussion of $\GAMExt$ and this last statement, referring the reader to \cite{GAMKrMe13} for the details.

\medskip
In order to discuss ellipticity of a boundary value problem we need to revisit the normal family, see \eqref{NormalFamilyAsLimit}. For each $y$, the operators $A_\wedge(\pmb\eta)$, $\pmb\eta\in T_y^*\GAMY$, are cone operators, elements of $x^{-1}\GAMDiff^1_b(\GAMZ^\wedge_y;E_{\GAMZ^\wedge_y},F_{\GAMZ^\wedge_y})$. As such they have maximal and minimal domains, and each its own trace space, which is just the fiber  $\GAMtrb_y$ lifted via $\pi_\GAMY:T^*\GAMY\to\GAMY$ to $\pmb\eta$ because the indicial operator of $A_\wedge(\pmb \eta)$ is $\GAMA_y$. They also have each a trace map $\gamma_{A_\wedge(\pmb\eta)}:\GAMDom_{\max}(A_\wedge(\pmb\eta))\to \pi_\GAMY^*\GAMtrb_y$. The vector bundle $\GAMKK\to T^*\GAMY\GAMminus 0$ whose fiber at $\pmb\eta$ is the kernel of
\begin{equation*}
A_{\wedge}(\pmb\eta) : \GAMDom_{\max}(A_{\wedge}(\pmb\eta))\subset x^{-1/2}L^2_b(\GAMZ^\wedge_y;E_{\GAMZ^\wedge_y}) \to x^{-1/2}L^2_b(\GAMZ^\wedge_y;F_{\GAMZ^\wedge_y}), y=\pi_\GAMY(\pmb\eta),
\end{equation*}
is a smooth vector bundle. The image $\gamma_{A_\wedge}(\GAMKK)=\GAMkerb\subset \pi^*_\GAMY\GAMtrb$ is a subbundle, see Theorem~6.2, {\GAMopcit}

\medskip
Boundary conditions will be, as generically in the classical case, conditions of the form $B\gamma_Au=g$ where
\begin{equation*}
B : C^{\infty}(\GAMY;\GAMtrb) \to C^{\infty}(\GAMY;\GAMtarget),\quad B\in \Psi^{\mu}_{1,\delta}(\GAMY;(\GAMtrb,-({x\partial_x+1/2})),(\GAMtarget,\GAMagen))
\end{equation*}
for some $\mu \in \GAMR$. As the notation indicates, the vector bundle $\GAMtarget\to\GAMY$ comes equipped with a smooth endomorphism $\GAMagen$. We assume that $B$ has twisted homogeneous principal symbol $\GAMsym(B)$ (see the end of Section~\ref{sec-variableorder} above). To account for the possible necessity of Atiyah-Patodi-Singer conditions, let 
\begin{equation*}
\Pi \in \Psi^0_{1,\delta}(\GAMY;(\GAMtarget,{\GAMagen}),(\GAMtarget,{\GAMagen}))
\end{equation*}
be a projection with twisted homogeneous principal symbol $\GAMsym(\Pi)$. On account that $\Pi$ is a continuous projection, its range, $\Pi H^{1-\mu+{\GAMagen}}(\GAMY;\GAMtarget)$, is a closed subspace of $H^{1-\mu+{\GAMagen}}(\GAMY;\GAMtarget)$. And because $\GAMsym(\Pi) = \GAMsym(\Pi)^2$, the range of $\GAMsym(\Pi)$ is a subbundle $\GAMtarget_{\Pi}$ of $\pi_{\GAMY}^*\GAMtarget$ over $T^*\GAMY\GAMminus 0$.

We are now ready to state the central result of \cite{GAMKrMe13}:

\begin{theorem}[Theorem 9.6, \GAMopcit]
Let $A\in x^{-1}\GAMDiff^1_e(\GAMM;E,F)$ be $w$-elliptic. Assume that \eqref{FiniteSpecb} holds with $\gamma=1/2$ and $m=1$, so there is a well defined trace bundle $\GAMtrb\to\GAMY$. Assume further that
\begin{equation*}
\begin{gathered}
A_{\wedge}(\pmb\eta) :\GAMDom_{\min}(A_{\wedge}(\pmb\eta)) \subset x^{-1/2}L^2_b(\GAMM;E) \to x^{-1/2}L^2_b(\GAMM;F) \text{ is injective}, \\
A_{\wedge}(\pmb\eta) :\GAMDom_{\max}(A_{\wedge}(\pmb\eta)) \subset x^{-1/2}L^2_b(\GAMM;E) \to x^{-1/2}L^2_b(\GAMM;F) \text{ is surjective}
\end{gathered}
\end{equation*}
for all $\pmb\eta \in T^*\GAMY\GAMminus 0$. Finally, assume
\begin{equation*}
\GAMsym(\Pi B) : \GAMkerb \to \GAMtarget_{\Pi}
 \text{ on } T^*\GAMY\GAMminus 0 \text{ is an isomorphism}.
\end{equation*}
Then the operator
\begin{equation*}
\begin{bmatrix} A \\ \Pi B \gamma_A \end{bmatrix} : H^1_{\GAMtrb}(\GAMM;E) \to
\begin{array}{c} x^{-1/2}L^2_b(\GAMM;F) \\ \oplus \\ \Pi H^{1-\mu+{\GAMagen}}(\GAMY;\GAMtarget) \end{array}
\end{equation*}
is a Fredholm operator.
\end{theorem}

In other words, the boundary value problem
\begin{equation*}
\left\{\begin{aligned}
&Au = f \in x^{-1/2}L^2_b(\GAMM;F)\\
&\Pi B(\gamma_A u) = g \in \Pi H^{1-\mu+{\GAMagen}}(\GAMY;\GAMtarget),
\end{aligned}
\right.
\end{equation*}
in which $u$ is sought in $H^1_{\GAMtrb}(\GAMM;E)$, is well-posed.

%%%%%%%%%%%%%%%%%%%%%%%%%%%%%%%%%%%%%%%%%%%%%%%%%

\end{document}